\documentclass[11pt,twoside,leqno,euscript]{aomamlt2e} 
  \pageno{715}
\received{September 16, 2004}

 \newcommand{\GL}{\mathrm{GL}}
\newcommand{\righthookdown}{\raise10pt%
  \hbox{$\scriptscriptstyle\cap$}%
  \kern-8.45pt\big\downarrow}
 \newcommand{\dbQ}{{\mathbb{Q}}}
\newcommand{\dbZ}{{\mathbb{Z}}}
\newcommand{\Kbar}{\bar{K}}
\newcommand{\kbar}{\bar{k}}
\newcommand{\Vbar}{\bar{V}}
\newcommand{\Ybar}{\bar{Y}}
\numberwithin{equation}{section}

\newcommand{\mathscr}{\EuScript}

\newcommand{\scrC}{{\mathscr{C}}} 
\newcommand{\scrL}{{\mathscr{L}}}

\newtheorem{thm}{Theorem}

\newtheorem{cor}[thm]{Corollary}

\theorembodyfont{\upshape}
\newtheorem{defn}[thm]{{\it Definition}}

\newtheorem{rmk}[thm]{{\it Remark}}

\numberwithin{thm}{section}

\newcommand{\ml}[2]{\begin{multline}\label{#1}#2 \end{multline}}
\newcommand{\ga}[2]{\begin{gather}\label{#1}#2 \end{gather}}

\newcommand{\surj}{\twoheadrightarrow}

\newcommand{\Spec}{{\rm Spec \,}}

\newcommand{\Gal}{{\rm Gal}}

\newcommand{\sC}{{\mathcal C}}

\newcommand{\sI}{{\mathcal I}}

\newcommand{\sL}{{\mathcal L}}
\newcommand{\sM}{{\mathcal M}}

\newcommand{\sO}{{\mathcal O}}

\newcommand{\F}{{\mathbb F}}

\newcommand{\N}{{\mathbb N}}
\renewcommand{\P}{{\mathbb P}}
\newcommand{\Q}{{\mathbb Q}}

\newcommand{\Z}{{\mathbb Z}}

\begin{document}
\currannalsline{164}{2006} 
 \newcommand{\ee}{{\hskip1pt\rm \'{\hskip-6.5pt \it e}}}
\newcommand{\gra}{{\hskip1pt\rm \`{ \hskip-6.5pt\it a}}}
 \newcommand{\gre}{{\hskip1pt\rm \`{ \hskip-6.5pt\it e}}}
\newcommand{\gu}{{\hskip1pt\rm \`{ \hskip-6.5pt\it u}}}
 \newcommand{\he}{{\hskip1pt\rm \^{\hskip-6.5pt \it e}}}
 \newcommand{\ho}{{\hskip1pt\rm \^{\hskip-6.5pt \it o}}}

 \title{Deligne's integrality theorem\\ in unequal characteristic and\\ rational points over finite fields}

 \acknowledgements{Partially supported by the 
DFG-Schwerpunkt ``Komplexe Mannigfaltigkeiten'' and by the DFG Leibniz Preis.}
 \author{H\'el\`ene Esnault}

 \institution{Universit\"at Duisburg-Essen, Mathematik, Essen, Germany\\
\email{esnault@uni-due.de}}

\shorttitle{Integrality}

\begin{center}{\it \`A Pierre Deligne{\rm ,} \gra\
 l\/{\rm '}\/occasion de son $60$-i\gre me anniversaire{\rm ,}\\
 en t\ee moignage de
profonde admiration}
\end{center}
 
\centerline{\bf Abstract}
\vglue8pt
If $V$ is a smooth projective variety defined over a local field $K$ 
with finite residue field, so that its \'etale cohomology over the algebraic closure  $\bar{K}$ 
is supported in codimension 1, 
then the mod $p$ reduction of a projective regular model
carries a rational point. As a consequence, if the Chow group of 0-cycles
of $V$ over a large algebraically closed field is trivial, then 
the mod $p$ reduction of a projective regular model
carries a rational point.

\section{Introduction}

If $Y$ is a  smooth, geometrically irreducible, projective
 variety over a finite field $k$, we singled out in \cite{Epoint} a motivic condition forcing the existence of a rational point. Indeed, if the Chow group of 0-cycles of $Y$ fulfills base change $CH_0(
Y\times_k \overline{k(Y)})\otimes_{\Z} \Q=\Q$, then the number of rational points of $Y$ is congruent to 1 modulo $|k|$.  In general it is hard to control the Chow group of 0-cycles, but if $Y$ is rationally connected, for example  if $Y$ is a 
Fano variety, 
then the base change condition is fulfilled, and thus, rationally connected varieties over a finite field have a rational point. Recall the Leitfaden of the proof. By S. Bloch's decomposition of the diagonal acting on cohomology as a correspondence 
\cite[Appendix to 
Lecture 1]{Bl}, the base change condition implies that \'etale cohomology $H^m(\bar{Y}, \Q_\ell)$ is
supported in codimension
$\ge 1$ for all $m\ge 1$, that is  that \'etale cohomology for $m\ge 1$ 
lives in 
coniveau 1. Here $\ell$ is a prime number not dividing $|k|$. On the other hand, by 
 Deligne's integrality theorem \cite[Cor.~5.5.3]{DeInt}, the coniveau condition implies that the
eigenvalues of the geometric Frobenius acting on $H^m(\bar{Y}, \Q_\ell)$ are divisible by  
\pagebreak $|k|$ as algebraic   integers for $m\ge 1$; thus the Grothendieck-Lefschetz  trace formula
\cite{Gr}  allows us 
to conclude.  Summarizing, we see that the cohomological condition which forces the existence of a rational
point is the coniveau condition. The motivic  condition  is here to allow us  to check geometrically in 
concrete examples the coniveau condition.

If $Y$ is no longer smooth, then homological
cycle classes no longer act on cohomology; thus the base change condition is no longer the right condition to force the existence of a rational point. Indeed, J. Koll\'ar constructed  an example of a 
 rationally connected projective variety, but without any rational 
 point. On the other hand, the classical theorem by Chevalley-Warning \cite{Ch}, \cite{War}, and its generalization by Ax-Katz \cite{Ax}, \cite{Ka}, asserting that the number of rational points of a closed subset $Y$ of $\P^n$ defined by $r$ equations of degree $d_i$, with $ \sum_1^r d_i\le n$, is congruent to 1 modulo $|k|$, suggests that when  $Y$ is smoothly deformable, the rational points of the smooth fibres singled out in \cite{Epoint}  produce rational points on the singular fibres of the deformation. Indeed, N. Fakhruddin and C. S. Rajan showed that if $f: X\to S$ is a projective dominant morphism over a finite field, with $X,S$ smooth connected, and if the base change condition is generically
 satisfied, that is if $CH_0(X\times_{k(S)} \overline{k(X)})
\otimes_{\Z}\Q=\Q$, then the number of rational points of a closed fibre is congruent to 1 modulo the cardinality of its field of definition
\cite[Th.~1.1]{FR}. The method is a refined 
version of the one explained above 
 in the smooth case, that is when $S$ is the spectrum of a finite field.
However,  it does not allow us to finish the proof  if only  
the coniveau condition on the geometric  
general fibre is known.  On the other hand, the previous discussion in the smooth case indicates that it should be sufficient 
to assume that the geometric general fibre fulfills the cohomological 
coniveau condition to force 
the singular fibres to acquire a rational point. According to 
Grothendieck's and Deligne's philosophy of motives, which links the  level for the congruence of rational points over finite fields to the level for the Hodge type over the complex numbers, 
this is supported by the fact that if $f: X\to S$ is a projective dominant
morphism 
over the field of complex numbers, with $X, S$ smooth, $S$ a connected curve, 
and if   the Hodge type of some smooth closed fibre is at least 1, 
then  so is the Hodge type of all closed fibres  \cite[Th.~1.1]{Eapp}).

We state now our theorem and several consequences.
 Let $K$ be a local field, with ring of integers $R\subset K$ and finite residue field $k$. 
We choose 
a prime number $\ell$
not dividing $|k|$.  If $V$ is a variety defined over $K$, 
we denote by $H^m(V\times_K \bar{K}, \Q_\ell)$ its $\ell$-adic cohomology. We say that $H^m(V\times_K \bar{K}, \Q_\ell)$
has coniveau 1 if each class in this group dies
in $H^m(U\times_K \bar{K}, \Q_\ell)$  after restriction 
on some nonempty open $U\subset V$. 

 \vglue-36pt
\phantom{up}
\begin{thm} \label{mainthm} Let $V$ be an absolutely irreducible{\rm ,} 
 smooth projective variety over $K${\rm ,}
 with a regular projective model $X$ over $R$. If \ee tale cohomology $H^m(V\times_K \bar{K}, \Q_\ell)
$ 
has coniveau $1$ for all $m\ge 1${\rm ,}
 then the number of rational points of the special fibre $Y=X\times_R k$ is congruent to $1$ modulo $|k|$.
\end{thm}

 Let $K_0 \subset K$ be a subfield of finite type over its prime field 
over which $V$ is defined, i.e. $V=V_0\times_{K_0}K$ for some variety $V_0$ defined over $K_0$, 
 and let $\Omega$ be a field extension  of $K_0(V_0)$. For example if $K$ has unequal characteristic, we may take $\Omega=K$. 
Using the decomposition of the diagonal mentioned before, one obtains  

\begin{cor} \label{mainthm:cor}
 Let $V$ be an absolutely irreducible{\rm ,}
 smooth projective variety over $K${\rm ,} with a regular projective model $X$ over $R$.
 If the Chow group of
$0$-cycles fulfills base change
$CH_0(V_0\times_{K_0} \bar{\Omega})\otimes_{\Z} \Q=\Q${\rm ,} then the number
of rational points of $Y$ is congruent to $1$ modulo $|k|$. \end{cor}

(See  \cite[Question 4.1]{FR} for the corollary, 
where the regularity of $X$ is not asked for.) 
In particular, our corollary applies for Fano varieties, and more generally, for rationally connected varieties $V$.

If the local field $K$ has equal characteristic, this is a certain strengthening of  
\cite[Th.~1.1]{FR}. Indeed,  our basis ${\rm Spec}(R)$
has only Krull dimension 1, but our coniveau assumption is the one which was expected, as indicated above.
If the local field 
$K$ has unequal characteristic, 
we see directly Deligne's philosophy at work. To our knowledge, this is the first such example. 
In this case,  the coniveau $ 1$ condition for \'etale cohomology is equivalent to the coniveau $ 1$ condition for de Rham cohomology $H^m_{\rm DR}(V\times_K \bar{K})$. 
By
Deligne's mixed Hodge theory \cite{DeHodgeII}, it implies that
 the Hodge type of de Rham cohomology $H^m_{\rm DR}(V)$ is $\ge 1$ for all $m\ge 1$, or equivalently that $H^m(V, \sO_V)=0$ for all $m\ge 1$. Conversely, Grothendieck's generalized conjecture predicts that 
those two conditions are equivalent; that is the Hodge type being $\ge 1$ should imply that the coniveau is  1. Thus one expects that if 
$V$ is a smooth projective variety over $K$, with 
$H^m(V, \sO_V)=0$ for all $m\ge 1$, then if $X$ is a regular projective model of $V$, 
the number of rational points of $Y=
X\times_R k$  is congruent to 1 modulo $|k|$. 
In particular this  holds for surfaces. 
\begin{thm} \label{surf:thm} Let $V$ be an absolutely irreducible{\rm ,}
 smooth projective
surface defined over a finitely generated $\Q$-algebra  $L$. If 
$$H^1(V, \sO_V)=H^2(V, \sO_V)=0,$$ then for any prime place
of $L$ with $p$-adic completion $K${\rm ,}
 for which $V\times_L K$ has a regular model $X${\rm ,}
 the number of rational points of the  mod $p$ reduction 
$X\times_R k${\rm ,}
 where $R\subset K$ is the ring of integers and $k$ is the finite residue field{\rm ,} is congruent to $1$  
modulo $|k|$. 
\end{thm} 

An example of such a surface is Mumford's fake $\P^2$ \cite{Mu}, a surface in 
characteristic 0 which has the topological invariants of $\P^2$, yet is of general type. We still do not know whether its Chow 
\pagebreak group of 0-cycles fulfills base change, as predicted by Bloch's conjecture. The surface is
constructed
 by\break 2-adic uniformization,
and the special fibre over $\F_2$, says Mumford quoting Lewis Carroll to 
express his ``confusion'', is a $\P^2$ blown up 7 times, crossing itself in 7 rational double curves, themselves crossing in 7 triple points \ldots Theorem \ref{surf:thm} allows one to say
 (in a less 
entertaining way) that at other bad primes with a regular projective model, there are rational points as well.

We now describe our method. Our goal is to show that the eigenvalues of the geometric Frobenius $F\in {\rm Gal}(\bar{k}/k)$ acting on 
$H^m(Y\times_k \bar{k},\Q_\ell)$ are $|k|$-divisible algebraic integers
for $m\ge 1$. 
Indeed, this will imply, by the Grothendieck-Lefschetz trace formula \cite{Gr},
that  $Y$ has modulo $|k|$ the same number
of  
rational points as $\P^N_k$.

 To this aim, we  consider the specialization map 
$H^m(Y\times_k \bar{k},\Q_\ell)
\xrightarrow{\rm sp} H^m(V\times_K\bar{K}, \Q_\ell)$ which is the edge 
homomorphism in the vanishing cycle spectral sequence (\cite[p.\ 214, (7)]{DeWeII}, 
\cite[p.\ 23]{RZ}). 
 Let $G$ be the Deligne-Weil group of the local field $K$. 
This is an extension of $\Z$, generated multiplicatively by the geometric Frobenius $F$ of ${\rm
Gal}(\bar{k}/k)$, by the inertia $I$. It acts on $H^m(V\times_K \bar{K}, \Q_\ell)$, on
$H^m(Y\times_k \bar{k},\Q_\ell)$ via its quotient $\Z\cdot F$, and the specialisation map is
$G$-equivariant. On the other hand, denoting by
$K^u$    the maximal unramified extension of $K$ in $\bar{K}$, that is 
$K^u=K^I$, the 
specialization map has a  $G$-equivariant factorization 
$$sp: H^m(Y\times_k \bar{k},\Q_\ell)
\xrightarrow{sp_u} H^m(V\times_K K^u, \Q_\ell)
\to  H^m(V\times_K\bar{K}, \Q_\ell),$$ where on the first two terms,
$G$ acts via its quotient $\Z\cdot F$. 
We first show 
\begin{thm} \label{thm:ker} Let $V$ be a smooth projective variety over a local field $K$ with finite residue field $k$ .
If $X$ is a regular projective model over $R${\rm ,} then 
the eigenvalues of $F$ on  
the kernel of the specialization map $sp_u$ 
are $|k|$-divisible  algebraic integers.
\end{thm}

Theorem \ref{thm:ker} is a consequence of Deligne's integrality theorem 
loc. cit. and of  Gabber's purity theorem \cite[Th.~2.1.1]{Fu}.

This reduces the problem to showing $|k|$-divisibility of the eigenvalues of 
$F$ on ${\rm Im}( sp_u) \subset 
H^m(V\times_K K^u, \Q_\ell) $. The latter group is an $F$-equivariant
extension of the inertia invariants $H^m(V\times_K \bar{K},\Q_\ell)^I$ by the first inertia cohomology group $H^1(I, H^{m-1}(V\times_K 
\bar{K},\Q_\ell))$. 
By Grothendieck \cite{Gr2}, as $k$ is finite, $I$ acts quasi-unipotently
on $H^m(V\times_K \bar{K}, \Q_\ell)$. As a consequence,
  modulo multiplication by roots of unity, 
the eigenvalues of a lifting $\Phi \in G$ of $F$ acting on
$H^m(V\times_K \bar{K},\Q_\ell)$
depend only on $F$ 
(\cite[Lemme (1.7.4)]{DeWeII}). In particular, if for one choice of $\Phi$, there are algebraic integers, then
they are algebraic integers for all choices.  We denote by $N^1H^m(V\times_K \bar{K}, \Q_\ell)$ the
subgroup of
$H^m(V\times_K \bar{K}, \Q_\ell)$ consisting of the classes which die
in $H^m(U\times_K \bar{K}, \Q_\ell)$
 after restriction on some nonempty open $U\subset V$.  It is a $G$-submodule. 
Then Theorem \ref{mainthm} is a consequence of 
\begin{thm} \label{thm:int}
Let  $V$ be a smooth irreducible projective variety defined over a local field $K$ with finite residue field $k$. Let $\Phi$ be a lifting of 
the geometric Frobenius of $k$ in the Deligne-Weil group of $K$. Then 
the eigenvalues of $\Phi$
\begin{itemize}
\item[{\rm i)}]   on $H^m(V\times_K \bar{K}, \Q_\ell)$ 
are algebraic integers for all $m${\rm ,}
\item[{\rm ii)}] on
 $N^1H^m(V\times_K 
\bar{K}, \Q_\ell)$  are 
$|k|$-divisible algebraic integers. \end{itemize} 
\end{thm}

Theorem \ref{thm:int} is a consequence of Deligne's integrality theorem 
loc.\ cit., of  de Jong's alterations \cite{deJ} and of Rapoport-Zink's weight spectral sequence \cite{RZ}.

\demo{Acknowledgements} We thank Pierre Berthelot, Gerd Faltings
for discussions, Jean-Louis Colliot-Th\'el\`ene and  
Wayne Raskind for careful reading of an earlier version of the article and for comments. 
We heartily thank Spencer Bloch for suggesting that we
  compute on $K^u$ and for his encouragement,\break  Johan de Jong
for pointing out an error in the proof of Theorem \ref{mainthm} in the  first version of the article, and the referee for  forcing  and helping us to restore the whole strength of Theorem \ref{mainthm} in the corrected version.
We thank the Alfr\'ed R\'enyi Institute, Budapest, for its support during the preparation of part of this work.

\section{The kernel of the specialization map over\\ the maximal unramified extension}

Let $V$ be a smooth projective variety over a local field $K$ with projective model $X$ over the ring of
integers and special fibre $Y=X\times_R k$ over the residue field $k$ which we assume throughout to be
finite.

In the following, $K^u$ is the maximal unramified extension of $K$, $R^u$ its ring of integers, with residue field $\bar{k}$.   
The specialization map $sp_u$ \cite[p.~213 (6)]{DeWeII},  is 
constructed by applying base change $H^m(Y\times_k \bar{k}, \Q_\ell)
=H^m(X\times_R R^u, \Q_\ell)$ for $X$ proper over $R$, followed by the restriction map
$ H^m(X\times_R R^u, \Q_\ell) \to  
H^m(V\times_K K^u, \Q_\ell)$.  In particular, one has an exact sequence
\begin{eqnarray}
\label{2.1}
 \ldots \to H^m_Y(X\times_R R^u, \Q_\ell) &\to &H^m(Y\times_k \bar{k}, \Q_\ell) 
\xrightarrow{sp_u} H^m(V\times_K K^u, \Q_\ell)\\
&\to& H^{m+1}_Y(X\times_R R^u, \Q_\ell)\to \ldots\  . \nonumber
\end{eqnarray}
Here in the notation: $H_Y(()\times_R R^u,())$ means $H_{Y\times_R R^u}(()\times _R R^u, ())$ etc. 
The geometric Frobenius $F\in {\rm Gal}(\bar{k}/k)$ acts
on all terms in \eqref{2.1} and the exact sequence is $F$-equivariant. 
Theorem \ref{thm:ker}  is then a  consequence of 
\begin{thm} \label{gabber:thm}
If $X$ is a regular scheme defined over $R${\rm ,}
 with special fibre $Y=X\times_R k${\rm ,} then the eigenvalues of $F$ acting on 
$H^m_Y(X\times_R R^u, \Q_\ell)$ are algebraic integers in $|k|\cdot \bar{\Z}$ for all $m$. 
\end{thm}

\Proof 
We proceed as in \cite[Lemma 2.1]{Epoint}. One has a finite stratification 
$\ldots \subset Y_i \subset Y_{i-1} \subset \ldots \subset Y_0=Y$ 
by closed subsets defined over $k$ such that $Y_{i-1}\setminus Y_i$ is smooth.
It yields    the $F$-equivariant localization  sequence 
\begin{eqnarray}
\label{2.2} 
\ldots \to H^m_{Y_i}(X\times_R R^u, \Q_\ell)& \to& 
 H^m_{Y_{i-1}}(X\times_R R^u, \Q_\ell) \\ &\to& H^m_{(Y_{i-1}\setminus Y_i)}(
(X\setminus Y_i)\times_R R^u, \Q_\ell)\to \ldots\ .
\nonumber
\end{eqnarray}
Thus Theorem \ref{gabber:thm} is a consequence of 
\begin{thm} \label{gabber2:thm} If $X$ is a regular scheme defined over $R${\rm ,} 
and $Z\subset Y=X\times_R k$ is a smooth closed subvariety defined over $k${\rm ,} 
then the eigenvalues of 
$F$ acting on $H^m_Z(X\times_R R^u, \Q_\ell)$
lie in $|k|\cdot \bar{\Z}$ for all $m$. \end{thm} 

\Proof  The scheme $X$ defined over $R$ being regular, its base change $X\times_R R^u$ by the unramified  map ${\rm Spec} \ R^u\to {\rm Spec} \ R$ is regular as well. 
By Gabber's purity theorem \cite[Th.~2.1.1]{Fu}, one has an 
$F$-equivariant isomorphism  
\ga{2.3}{H^{m}(Z\times_k \bar{k}, \Q_\ell)(-c)\cong 
H^{m+2c}_Z(X\times_R R^u, \Q_\ell),
} 
where $c$ is the codimension of $Z$ in $X$. 
Thus in particular, $F$ acts on\break $H^{m+2c}_Z(X\times_R R^u, \Q_\ell)$  
 as it does on $H^{m}(Z\times_k \bar{k}, \Q_\ell)(-c)$. We are back to a problem over finite fields. Since $c\ge 1$, we only need to know that the eigenvalues of $F$ on $H^{m}(Z\times_k \bar{k}, \Q_\ell)$ lie in $\bar{\Z}$. This is 
\cite[Lemma 5.5.3 iii]{DeInt} (via duality as $Z$ is smooth).
\Endproof\vskip4pt  
This finishes the proof of Theorem \ref{gabber:thm}. 
\hfill\qed 

\begin{rmk} \label{rmk}
We observe that \eqref{2.1} together with Theorem \ref{gabber:thm} implies that if $V$ is a smooth projective variety defined over a local field $K$, and $V$ admits a regular model over $R$, then the eigenvalues of $F$ on
$H^m(V\times_K K^u, \Q_\ell)$ are algebraic integers, and they are $|k|$-divisible algebraic integers for
some $m$ if and only if the eigenvalues of $F$ on $H^m(Y\times_k \bar{k}, \Q_\ell)$ are $|k|$-divisible
algebraic integers
 for the same  $m$. 
\end{rmk}

\section{Eigenvalues of a lifting of Frobenius on \'etale cohomology of smooth projective varieties}

Let $V$ be a smooth projective variety over a local field $K$ with projective model $X$ over the ring of integers $R$ and special fibre $Y=X\times_R k$ 
over the finite residue field $k$. Let $\Phi$ be a
lifting of Frobenius in the Deligne-Weil group of $K$. The aim of this section is to prove Theorem \ref{thm:int}.

Recall that $X/R$ is said to be strictly semi-stable
 if $Y$ is reduced and is a strict normal crossing divisor. In this case, $X$ is necessarily regular as well.
Recall from \cite[(6.3)]{deJ}  that if $A \subset X, A=\sum_i A_i$ is a divisor, $(X,A)$ is said to be a strictly
semi-stable pair if 
$X/R$ is strictly semi-stable,
$A+Y$ is a normal crossing divisor, and all the strata $A_I/R, I=(i_1,\ldots, i_s)$ of
$A$ are strictly semi-stable as well. 

\demo{Proof of Theorem {\rm \ref{thm:int} i)}} \hskip-9pt Let $V$ be as in Theorem 
\ref{thm:int} i). Let \hbox{$K'\!\supset\! K$} be a finite extension, 
with residue field $k'\!\supset\! k$, and Deligne-Weil group \hbox{$G'\!\subset\! G$.}\break 
Let $\sigma: V'\to V$ be an alteration; that is, $V'$ is smooth projective over $K'$, $\sigma$ is proper, dominant  and generically finite. Then 
$\sigma^*: H^m(V\times_K \bar{K}, \Q_\ell)\to 
H^m(V'\times_{K'} \overline{K'}, \Q_\ell)$ is injective, and $G'$-equivariant. In particular, it\break is
$\Phi'$-equivariant for a lifting $\Phi'\in G'$ of $F^{[k':k]}$. Thus Theorem \ref{thm:int} for $\Phi'$
implies  Theorem \ref{thm:int} for $\Phi$. By de Jong's fundamental alteration theorem\break
(\cite[Th.~6.5]{deJ}), we may find such $K', V'$ with the property that $V'$ has a strictly semi-stable model 
over the ring of integers of $K'$.  Thus by the above, without loss of generality, we may assume that $V$
defined over $K$ has a strictly semi-stable model $X$ over the ring of integers $R\subset K$. We denote by
$Y=X\times_Rk$ the closed fibre. It is a strict normal crossing divisor. We denote by $Y^{(i)}$ the disjoint
union of the smooth strata of codimension $i$ in $X$. Thus $Y^{(0)}=X$, $Y^{(1)}$ is the disjoint union of
the  components of $Y$ etc.     We apply now the existence of the weight spectral sequence \cite[Satz
2.10]{RZ} by Rapoport-Zink  (see also \cite[(3.6.11), (3.6.12)]{I} for a r\'esum\'e),
\ml{3.1}{ 
{}_WE_1^{-r, m+r}=
\oplus_{q\ge 0, r+q\ge 0} H^{m-r-2q}(Y^{(r+1+2q)}
\times_k \bar{k}, \Q_\ell)(-r-q)\\
 \Rightarrow H^m(V\times_K \bar{K}, \Q_\ell)\ .
}
It is $G$-equivariant and converges in $E_2$ (\cite[p.~41]{I}).
Thus  the eigenvalues of $\Phi$ on the right-hand side are (some of) the eigenvalues of $F$ on the left-hand side.  We apply again Deligne's integrality theorem \cite{DeInt},
loc.\ cit.\  to conclude the proof. 
\Endproof\vskip4pt  

{\it Proof   of Theorem} \ref{thm:int} ii).
Let $V$ be as in Theorem \ref{thm:int} ii). 
Since \'etale cohomology  $H^m(V\times_K \bar{K}, \Q_\ell)$ is a finite
dimensional $\Q_\ell$-vectorspace, there is a  divisor $A^0$ defined over $K$ 
with a $G$-equivariant surjection
$H^m_{A^0}(V\times_K \bar{K}, \Q_\ell) \surj N^1H^m(V\times_K \bar{K}, \Q_\ell).$ 
  Let $K'\supset K$ be a finite extension, let $\sigma: V'\to V$ be an alteration. Then 
\ml{3.2}{
\sigma^*({\rm Im} (H^m_{A^0}(
V\times_K \bar{K}, \Q_\ell)) \\
\subset {\rm Im} (H^m_{\sigma^{-1}(A^0)}(
V'\times_{K'} \overline{K'}, \Q_\ell))\subset H^m(V'\times_{K'} \overline{K'}, \Q_\ell).
}
Since as in the proof of i), $\sigma^*: H^m(V\times_{K} \bar{K}, \Q_\ell)
\to  H^m(V\times_{K'} \overline{K'}, \Q_\ell)$ is $G'$-equivariant and injective, Theorem \ref{thm:int} ii) for 
$\Phi'$ implies Theorem 
\ref{thm:int} ii) for $\Phi$. 
We use again  de Jong's  alteration theorem \cite[Th.~6.5]{deJ}.   
There is a
finite extension $K'\supset K$, with an alteration $\sigma: V'\to V$ such that
$V'$ has a strict semi-stable model $X'$ over $R'$, the ring of integers of 
$K'$,  and is such that  the Zariski closure $A'$ of $\sigma^{-1}(A^0)$ in $X'$ has the property that $(X',A')$ is strictly semi-stable. 
Thus by the above, we may assume that $(X,A)$ is a strictly semi-stable pair, 
where $A$ is the Zariski closure of $A^0$ in $X$. If $I$ is a sequence 
$(i_1,i_2,\ldots, i_a)$ of pairwise distinct indices, we denote by $A_I$ the intersection $A_{i_1}\cap A_{i_2}\cap \ldots \cap A_{i_a}$.
One has the $G$-equivariant Mayer-Vietoris spectral sequence 
\ga{3.3}{E_1^{-a+1,b}=
\oplus_{|I|=a} H^{b}_{A_I}(V\times_K \bar{K}, \Q_\ell)
\Rightarrow H^{1-a+b}_A(V\times_K \bar{K}, \Q_\ell)}
together with the $G$-equivariant purity isomorphism (e.g.\ \cite[Th.~2.1.1]{Fu})
\ga{3.4}{ H^{b-2c_I}(A_I\times_R \bar{K}, \Q_\ell)(-c_I)\cong 
 H^b_{A_I}(V\times_K \bar{K}, \Q_\ell),}
where $c_I$ is the codimension of $A_I$ in $X$. 
Since $c_I\ge 1$, we conclude with Theorem \ref{thm:int} i). 
\hfill\qed

\section{The proof of Theorem \ref{mainthm} and its consequences}

{\it Proof of Theorem} \ref{mainthm}. 
We denote by $\Phi$ a lifting of Frobenius in the Deligne-Weil group of $K$. 
By  Remark \ref{rmk}, $|k|$-divisibility of the eigenvalues of $F$ acting on $H^m(Y\times_k \bar{k}, \Q_\ell)$ is equivalent to $|k|$-divisibility of the
eigenvalues of $\Phi$ acting on $H^m(V\times_k K^u, \Q_\ell)$. On the other hand, one has the $F$-equivariant 
exact sequence \cite[p.~213, (5)]{DeWeII}
\ml{4.1}{
0 \to H^{m-1}(V\times_K \bar{K}, \Q_\ell)_I(-1) \to 
H^m(V\times_K K^u, \Q_\ell) \\
 \to H^{m}(V\times_K \bar{K}, \Q_\ell)^I \to 0\ .
}
Here $_I$ means the inertia  coinvariants while $^I$ means the 
inertia invariants. The quotient map $H^{m-1}(V\times_K \bar{K}, \Q_\ell)
\surj 
H^{m-1}(V\times_K \bar{K}, \Q_\ell)_I$ is $\Phi$-equivariant.  Thus by 
Theorem \ref{thm:int} i), the eigenvalues of $F$ acting on 
$H^{m-1}(V\times_K \bar{K}, \Q_\ell)_I$ are algebraic integers for 
all $m$; thus on $H^{m-1}(V\times_K \bar{K}, \Q_\ell)_I(-1)$ they are $|k|$-divisible algebraic integers for all $m$. The injection $
 H^{m}(V\times_K \bar{K}, \Q_\ell)^I \hookrightarrow H^{m}(V\times_K \bar{K}, \Q_\ell)$ is $\Phi$-equivariant; thus by 
the coniveau assumption of Theorem \ref{mainthm} and 
Theorem \ref{thm:int} ii), the eigenvalues of $F$ acting on
$H^{m}(V\times_K \bar{K}, \Q_\ell)^I$ are $|k|$-divisible algebraic integers
for $m\ge 1$. 
Thus we conclude that the eigenvalues of $F$ acting on $H^m(Y\times_k \bar{k}, 
\Q_\ell)$ are $|k|$-divisible algebraic integers for all $m \ge 1$. 
 By the Grothendieck-Lefschetz trace formula \cite{Gr} applied to $Y$, this  shows that the number of rational points of $Y$ is congruent to 1
modulo~$|k|$. This finishes the proof of Theorem \ref{mainthm}.
\Endproof\vskip4pt  

{\it Proof of Corollary} \ref{mainthm:cor}.
One applies Bloch's decomposition of the diagonal \cite[Appendix to Lecture 1]{Bl}, as mentioned in the
introduction and detailed in \cite{Epoint}, in order to show that the base change condition on the Chow
group of 0-cycles implies the coniveau condition of Theorem \ref{mainthm}.  Indeed,
$CH_0(V_0\times_{K_0} \bar{\Omega})\otimes_{\Z}\Q=\Q$ implies  the existence of a  decomposition
$N\Delta\equiv \xi \times V_0 + \Gamma$  in $CH_{{\rm dim} (V)} (V_0\times_{K_0} V_0)$, where
$N\ge 1, N\in \N$, 
$\xi$ is a 0-cycle of $V_0$ defined over $K_0$, 
$\Gamma$ is a ${\rm dim}(V)$-cycle lying in $V_0\times_{K_0} A$, where $A$ is a divisor in $V_0$. This decomposition yields {\it a fortiori\/} a decomposition in 
$CH_{{\rm dim} (V)} ( (V\times_K V)\times_K \bar{K})$. 
The correspondence with $\Gamma$ has image in ${\rm Im} 
(H^m_A(V\times_K \bar{K}, \Q_\ell))\subset 
H^m(V\times_K \bar{K}, \Q_\ell)$, 
while the correspondence with $\xi\times V_0$ 
 kills $H^m(V\times_K \bar{K}, \Q_\ell)$  for $m\ge 1$
as it factors through the restriction
to $H^m(\xi\times_{K_0} \bar{K}, \Q_\ell)$. 
Thus $$N^1H^m(V\times_k \bar{K},
\Q_\ell)  =H^m(V\times_K \bar{K}, \Q_\ell).$$
We apply Theorem \ref{mainthm} to conclude the proof. 
\Endproof\vskip4pt  

{\it Proof of Theorem} \ref{surf:thm}.
In order to apply Theorem \ref{mainthm}, we just have to know that $H^1(V, \sO_V)=0$ is equivalent to the vanishing of de Rham cohomology $H^1_{\rm DR}(V)$. Thus by the comparison theorem, 
this implies 
$H^1_{{\rm \acute{e}t}}(V\times_K \bar{K}, \Q_\ell)=0$. Furthermore, $H^2(V, \sO_V)=0$ is equivalent to $N^1H^2_{\rm DR}(V)=H^2_{\rm DR}(V)$; 
thus by the comparison theorem, 
$
N^1H^2_{\rm \acute{e}t}(V\times_K \bar{K}, \Q_\ell)=
H^2_{\rm \acute{e}t}(V\times_K \bar{K}, \Q_\ell)$. Thus we can apply Theorem
\ref{mainthm}.
\hfill\qed
\section{Some comments and remarks}

5.1.  
Theorem \ref{thm:int} ii) is formulated for $N^1$ and not for the higher coniveau levels $N^\kappa$ of \'etale cohomology. The appendix to this article fills  in this gap: if $V$ is smooth over a local field $K$ with finite 
residue field $k$, then the eigenvalues of $\Phi$ on $N^\kappa 
H^m_{{\rm prim}}(V\times_K \bar{K}, \Q_\ell)$ lie in $|k|^\kappa \cdot \bar{\Z}$. Here the subscript prim means one mods out by the powers of the 
class of the polarization coming from a projective embedding $Y\subset \P^N$. So for example, in the good reduction case, the $N^\kappa$ condition on the smooth projective fibre $V$ will imply that $|Y(k)|\equiv |\P^N(k)|$ modulo $|k|^\kappa$. In general, only a strong minimality condition on the model $X$ could imply this conclusion, as blowing up a smooth point of $Y$ keeps the same number of rational points only modulo $|k|$.

\demo{{\rm 5.2}} Koll\'ar's example of a rationally connected surface (personal
communication) 
over a finite field $k$, but without a rational point, is birational (over~$\bar{k}$) to the product of a genus
$\ge 2$ curve with $\P^1$. In particular it is not a Fano variety. Here we define a projective variety $Y$ over
a field $k$ to be Fano if it is geometrically irreducible, Gorenstein, and if the dualizing sheaf $\omega_Y$ is
anti-ample.  If the characteristic of $k$ is 0, then one defines the  ideal sheaf $\sI=\pi_*\omega_{Y'/Y}$,
where $\pi: Y'\to Y$ is a  desingularization. This ideal does not depend on the choice of $Y'$ (and 
is called in our days the multiplier ideal). The Kawamata-Viehweg vanishing theorem applied to $\pi^*\omega_Y^{-1}$ shows that $H^m(Y, \sI)=0$, for all $m$ if
$\sI$ is not equal to $\sO_Y$, otherwise for $m\ge 1$. In the cases where the support $S$ of $\sI$ is the empty set
 or where $S$ equals the singular locus of $X$, this implies by
\cite[Prop.~1.2]{E}
that the Hodge type of $H^m_{\rm DR}(X,S)$ is $\ge 1$ for all $m$ 
if $\sI$ is not equal to $\sO_Y$, otherwise for $m\ge 1$.
Using again Deligne's philosophy as mentioned in the introduction, one would expect that a suitable definition of $S$ in positive characteristic for a Fano variety (note the definition above requires the existence of a desingularization) would lead to the prediction that over a finite field $k$, the number of rational points of $Y$ is congruent to the number of rational points of $S$ modulo $|k|$ if $S=\emptyset$ or $S$ equals the singular locus of $X$.

\demo{{\rm 5.3}}
The correct motivic condition for a projective variety defined $Y$ over a finite field $k$, which implies that the number of rational points of $Y$ is 
congruent to 1 modulo $|k|$, is worked out in \cite{BEL}. Let us briefly recall
  it. If $Y\subset \P^N$ is any projective embedding defined over $k$,
 with complement $U$, and $\sigma: \P\to \P^N$ is an alteration of $\P^N$ which makes
$\sigma^{-1}Y=:Z$ a normal crossing divisor, then there is a divisor $A\subset \P$ in good position
with respect to all strata of $Z$ so that rationally, the  motivic  class $\Gamma_\sigma
\in H^{2N}((\P\times U) \times_k \bar{k}, (Z\times U)\times_k \bar{k}, N)$ dies on $\P\setminus A$
\begin{eqnarray}
\label{5.1}
0&=&\Gamma_\sigma|_{\P\setminus A} \\
&\in& H^{2N}(((\P\setminus A)\times U) \times_k \bar{k}, ((Z \setminus Z\cdot A)\times U)\times_k
\bar{k}, N)\otimes_{\Z}\Q.\nonumber
\end{eqnarray}
 Here $\Gamma_\sigma$ is the graph of the alteration $\sigma$. And
we know by the main theorem in \cite{BEL} that hypersurfaces $Y\subset \P^N$ of degree $d\le N$ fulfill
this motivic condition. Given Corollary \ref{mainthm:cor}, it is conceivable that the condition
$CH_0(V_0\times_{K_0}
\bar{\Omega}) \otimes_{\Z}
\Q\break =\Q$ directly implies our motivic condition \eqref{5.1}. We could not compute this, as moving
cycles over discrete valuation rings is not easy.

\demo{{\rm 5.4}} In   light of \cite[Th.~1.1]{Eapp}, 
one would expect that if $V$ is projective smooth over $K$ with smooth model $X$ over $R$,  if $K$ has unequal characteristic, and if 
$H^m(V, \sO_V)=0$ for all $m\ge 1$, then 
there is a way using rigid geometry which allows us to conclude that $Y$ 
has a point. 
 Grothendieck's generalized Hodge conjecture surely predicts that the coniveau 1 condition on \'etale cohomology of $V\times_K \bar{K}$ and the vanishing of $H^m(V, \sO_V)$ for $m\ge 1$ are equivalent, but we do not know the answer so far. 
It would give some hope to link  higher Hodge levels  $\kappa$ for $H^m_{\rm DR}(V)$ to higher levels
$\kappa$ for divisibility of Frobenius eigenvalues in \'etale cohomology and to higher levels $\kappa$ for congruences for the number of points of $Y$. 
In particular, it would give a natural explanation of the main results in \cite{E2}
and \cite{EK} where it is shown that for a closed subset $Y\subset \P^N$ defined over a finite field, the divisibility of the eigenvalues of $F$ is controlled by the divisibility for the number of rational points of $\P^N\setminus Y$ as stated in the Ax-Katz theorem \cite{Ax}, \cite{Ka}. 

\demo{{\rm 5.5}}
 We give a concrete nontrivial  example of Theorem \ref{mainthm} due to X. Sun (personal
communication).  Moduli $M(C, r,L)$ 
of vector bundles of rank $r$ and fixed determinant $L$ of degree $d$ with $(r,d)=1$ on a smooth projective curve $C$ over a field are known to be smooth projective Fano varieties, to which we can apply our Theorem \cite{Epoint} if the field is finite. If $(C,L)$ is defined over the local field $K$, with model $(\sC, \sL)$ over $R$ and reduction $(C_k, L_k)$ over $k$, then if $C_k$ has a node and $d=1$, $ M(C,2,L)$ has a model $\sM(\sC,2, \sL)$ with closed
fibre $M(C_k, 2, L_k)$ such that the underlying reduced variety parametrizes torsion-free sheaves $E$ of rank 2 which are endowed with a morphism $\Lambda^2E\to L_k$ which is an isomorphism off the double point. By \cite{Epoint},    
there is a rational point $E_0\in M(\tilde{C}_k, 2, \tilde{L}_k)(k)$, where $
\pi: \tilde{C}_k\to C_k$ is the normalization,  $\tilde{L}_k=\pi^*(
L_k\otimes \frak{m})$, and ${\rm Spec}(\sO_{C_k}/\frak{m})$ is the node.  
Then $E=\pi_*E_0 \in M(C_k,2,L_k)(k) $ is the wanted rational point of the modulo $p$ reduction. 
 
{
}
 
\centerline{\small (Received May 27, 2004)}

 \font\thinlinefont=cmr5
\font\sans=cmss12
  
 \vglue18pt
\begin{center}
{\bf \Large Appendix}
  \end{center}
\setcounter{section}{0}

  \vglue5pt \centerline{\small By {\scshape Pierre Deligne} and {\scshape H\'el\`ene Esnault}} 
\renewcommand{\institution}[1]
{\renewcommand{\theinstitutions}{\vskip16pt\baselineskip10pt\begin{quote}
\scriptsize\scshape #1\end{quote}}}

 \institution{The Institute of Advanced Study, School of  Mathematics,  Princeton, NJ\\
{deligne@math.ias.edu}
\\
\vglue-9pt
Universit\"at Duisburg-Essen, Mathematik, Essen, Germany\\
\email{esnault@uni-due.de}
}

 \shortname{Pierre Deligne and H\'el\`ene Esnault}

\vglue12pt

 We generalize in this appendix Theorem 1.5 to nontrivial coefficients on varieties $V$ which are neither smooth nor projective. 
 We thank Alexander Beilinson,  Luc Illusie and Takeshi Saito
for very helpful discussions.

The notation is  as in the article. Thus $K$ is a local field with finite residue field $k$, $R\subset K$ is the ring of integers, $\Phi$ is a lifting of the geometric Frobenius in the Galois group of $K$. We consider
  $\ell$-adic sheaves on schemes of finite type defined over $K$ in the sense of 
\cite[(1.1)]{apDeWeII}. One generalizes the definition \cite[ D\'ef.~5.1]{apDeInt}
of $T$-integral $\ell$-adic  sheaves on schemes of finite type
 defined over finite fields to $\ell$-adic sheaves  on schemes of finite
type  defined over local fields with finite residue field. Recall, to this aim, that if $\scrC$ is an $\ell$-adic sheaf on a $K$-scheme $V$ 
and $v$ is a closed point of $V$, then 
the stalk $\scrC_{\bar{v}}$ of $\scrC$ at $\bar{v}$ is a ${\rm Gal}(\bar{K}/K_v)$-module, where $K_v\supset K$ is the residue field of $v$,
 with residue field $\kappa(v)\supset k$.   On $\scrC_{\bar{v}}$ the inertia $I_v={\rm Ker} ({\rm Gal}(\bar{K}/K_v)\to {\rm Gal}(\kappa(v)/k))$  acts quasi-unipotently (\cite{apGr2}).
 Consequently the eigenvalues
of a lifting $\Phi_v\in {\rm Gal}(\bar{K}/K_v) $ of the geometric Frobenius $F_v \in {\rm Gal}
(\overline{\kappa(v)}/\kappa(v))$ are, 
up to multiplication by roots of unity, well defined (\cite[Lemma (1.7.4)]{apDeWeII}). 
 Let $T\subset \Z$ be a set  of  prime numbers.
\begin{defn} \label{apdefn:int} The $\ell$-adic sheaf $\scrC$ is $T$-integral if the eigenvalues of $\Phi_v$ acting on $\scrC_{\bar{v}}$
are integral over $\Z[\frac{1}{t}, t\in T]$ for all closed points $v\in V$.
\end{defn}

\begin{thm} \label{apthm:int}
Let $V$ be a scheme of finite type  defined over  $K${\rm ,} and let $\scrC$ be a $T$-integral $\ell$-adic sheaf on $V$. Then if $f: V\to W$ is a morphism to another $K$-scheme of finite type  $W$ defined over $K${\rm ,} the 
$\ell$-adic  sheaves 
$R^if_!\scrC$ are $T$-integral as well.
More precisely{\rm ,} if $w\in W$ is a closed point{\rm ,} then the eigenvalues of both 
$F_w$ and $|\kappa(w)|^{n-i}F_w$ acting on $(R^if_!\scrC)_{\bar{w}}$ 
  are integral over $\Z[\frac{1}{t}, t\in T]${\rm ,} with $n={\rm dim}(f^{-1}(w))$.
\end{thm}
\Proof 
Let $w$ be a closed point of $W$.
By base change for $Rf_!$ and by
$\{w\}\hookrightarrow W$, one is reduced to the case
where $W$ is the spectrum of a finite extension $K'$
of $K$.
If $\Vbar:=V\otimes_{K'}\Kbar'$, for $\Kbar'$ an
algebraic closure of $K'$, one has to check the
integrality statements for the eigenvalues of a
lifting of Frobenius on $H_c^i(\Vbar,\scrC)$.

Let us perform the same reductions as in \cite[p.~24]{apDeInt}.
Note that in loc. cit. $U$ can be shrunk so as to be
affine, with the sheaf smooth on it ($=$ a local
system).
This reduces us to the cases where $V$ is of
dimension zero or is an affine irreducible curve
smooth over $W=\Spec(K')$ and $\scrC$ is a smooth 
sheaf.

The integrality statement to be proven is insensitive
to a finite extension of scalars $K''/K'$.
The $0$-dimensional case reduces in this way to the
trivial case where $V$ is a sum of copies of
$\Spec(K')$.
In the affine curve case, $H_c^0$ vanishes, while, as
in \cite[Lemma 5.2.1]{apDeInt}, there is a $0$-dimensional
$Z\subset V$ such that the natural
($\Phi$-equivariant) map from $H^0(Z,\scrC)(-1)$ to
$H_c^2(\bar{V},\scrC)$ is surjective, leaving us
only $H_c^1$ to consider.

Let $\scrC_{\dbZ_\ell}$ be a smooth $\dbZ_\ell$-sheaf
from which $\scrC$ is deduced by $\otimes\dbQ_\ell$,
and let $\scrC_\ell$ be the reduction modulo $\ell$
of $\scrC_{\dbZ_\ell}$.
For some $r$, it is locally (for the \'etale topology)
isomorphic to $(\dbZ/\ell)^r$.
Let $\pi\colon\, V'\to V$ be the \'etale covering of
$V$ representing the isomorphisms of $\scrC_\ell$
with $(\dbZ/\ell)^r$.
It is a $\GL(r,\dbZ/\ell)$-torsor over~$V$.
As $H_c^*(\Vbar,\scrC)$ injects into
$H_c^*(\overline{V'},\pi^*\scrC)$, renaming irreducible
components of $V'$ as $V$ and $K'$ as $K$, we may and
shall assume that $W=\Spec(K)$ and that $\scrC_\ell$
is a constant sheaf.

Let $V_1$ be the projective and smooth completion of
$V$, and $Z:=V_1\setminus V$.
Extending scalars, we may and shall assume that $Z$
consists of rational points and that $V_1$, marked
with those points, has semi-stable reduction.
It hence is the general fiber of $X$ regular and
proper over $\Spec(R)$, smooth over
$\Spec(R)$ except for quadratic nondegenerate
singular points, with $Z$ defined by disjoint
sections $z_\alpha$ through the smooth locus.

Let $Y$ be the special fiber of $X$, and $\Ybar$ be
$Y\times_k\kbar$, for $\kbar$ the residue field of
the algebraic closure $\Kbar$ of $K$.

$$
\begin{matrix}
V &\null\kern5.0 true cm
\hidewidth:\text{\small{complement of 
     disjoint sections $z_\alpha$}}\hidewidth &&&\\
\raise3pt\hbox{$\scriptstyle j$}
   \righthookdown &&&&\\
V_1 &\hookrightarrow &X &\hookleftarrow &Y\\
\big\downarrow 
&&\big\downarrow &&\big\downarrow\\
{\text{\small Spec}}(K) &\hookrightarrow
&{\text{\small Spec}}(R) &\hookleftarrow 
  &{\text{\small Spec}}(k)\ . \end{matrix}
$$
The cohomology with compact support
$H_c^1(\Vbar,\scrC)$ is  $H^1(\overline{V_1},j_!\scrC)$, and
vanishing cycles theory relates this $H^1$ to the
cohomology groups on $\Ybar$ of the nearby cycle
sheaves $\psi^i(j_!\scrC)$, which are $\ell$-adic
sheaves on $\Ybar$, with an action of $\Gal(\Kbar/K)$
compatible with the action of $\Gal(\Kbar/K)$
(through $\Gal(\kbar/k)$) on $\Ybar$.
The choice of a lifting of Frobenius, i.e. of a
lifting of $\Gal(\kbar/k)$ 
in $\Gal(\Kbar/K)$, makes them come
from $\ell$-adic sheaves on $Y$, to which the
integrality results of \cite{apDeInt} apply.
Using the exact sequence
$$
0\to H^1(\Ybar,\psi^0(j_!\scrC))\to
H^1(\overline{V_1},j_!\scrC)\to H^0(\Ybar,\psi^1(j_!\scrC))
$$
and \cite{apDeInt} Th\'eor\`eme 5.2.2, we are reduced to check integrality of
the sheaves $\psi^i(j_!\scrC)$ ($i=0,1$).
It even suffices to check it at any $k$-point $y$ of
$Y$, provided we do so after any unramified finite
extension of $K$.

Let $X_{(y)}$ be the henselization of $X$ at $y$, and
$Y_{(y)}$, $V_{1(y)}$ and $V_{(y)}$ be the inverse
image of $Y$, $V_1$ or $V$ in $X_{(y)}$.
There are three cases:
$$
\epsfig{file=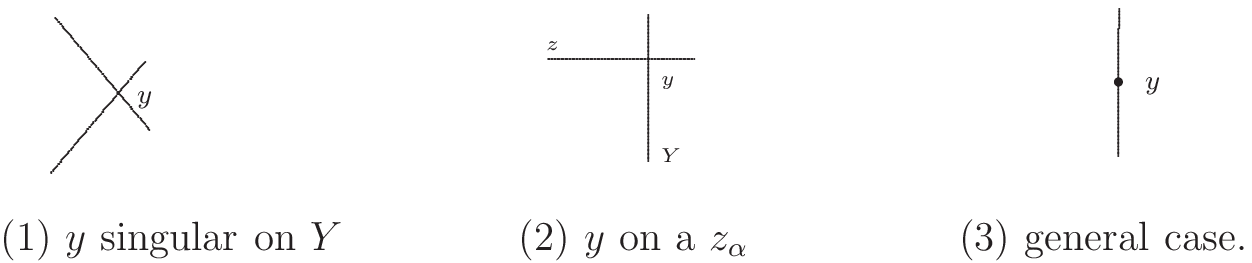}
$$
The restriction of $\psi^i$ to $y$ depends only on
the restriction of $\scrC$ to $V_{(y)}$, and short 
exact sequences of sheaves give rise to long exact
sequences of $\psi$.

Because $\scrC_\ell$ is a constant sheaf, $\scrC$ is
tamely ramified along $Y$ and the $z_\alpha$. 
More precisely, it is
given by a representation of the pro-$\ell$
fundamental group of $V_{(y)}$.
It is easier to describe the group deduced from the
profinite fundamental group by pro-$\ell$ completing
only the kernel of its map to
$\widehat{\dbZ}=\Gal(\kbar/k)$.
By Abhyankhar's lemma, this group is an extension of
$\widehat{\dbZ}$, generated by Frobenius, by
$\dbZ_\ell(1)^2$ in case (1) or (2) or $\dbZ_\ell(1)$
in case (3).
The representation is given by $r\times r$ matrices
congruent to $1\, {\rm mod}\, \ell$.
For $\ell\not=2$, such a matrix, if quasi-unipotent,
is unipotent.
Indeed, it is the exponential of its logarithm and
the eigenvalues of its logarithm are all zero.
For $\ell=2$, the same holds if the congruence is
${\rm mod}\,4$, hence if  $\scrC_{\dbZ_2}\, {\rm mod}\, 4$ is constant,
a case to which one reduces by the same argument we
used ${\rm mod}\, 2$.
By Grothendieck's argument\break \cite[p.~515]{apSeTa}, the action of
$\dbZ_\ell(1)$ or $\dbZ_\ell(1)^2$ is
quasi-unipotent, hence unipotent, and we can filter
$\scrC$ on $V_{(y)}$ by smooth sheaves such that
the successive quotients $Q$ extend to 
smooth sheaves on $X_{(y)}$.
If $Q$ extends to a smooth sheaf $\scrL$ on
$X_{(y)}$, the corresponding $\psi$ are known by
Picard-Lefschetz theory: $\psi^0$ is $\scrL$ restricted to
$Y_{(y)}$ in cases (1) and (3), and $\scrL$ outside
of $y$ extended by zero in case (2); $\psi^1$ is 
nonzero only in case (1), where it is $\scrL(-1)$ on
$\{y\}$ extended by zero.
By d\'evissage, this gives the required integrality.
\hfill\qed
\begin{cor}\label{apcor:int} Let $V$ be a scheme of finite type defined 
over $K$.  Then the eigenvalues of $\Phi$ on $H^i(\bar{V}, \Q_\ell)$ are integral over $\Z$. 
\end{cor}
\Proof  We fix an integer $n>2i$.
By de Jong's theorem \cite[Th.~6.5]{apdeJ}, there is a finite extension $K'\supset K$ and a simplicial truncated
alteration $V_n\to \ldots \to V_0\to V$ defined over $K'$. By the Mayer-Vietoris spectral sequence, $
H^*(\bar{V}, \Q_\ell)$ is filtered by sub ${\rm Gal}(\bar{K}/K')$-modules, the graduation of which is a
subquotient of $\oplus_0^n H^*(\bar{V}_j, \Q_\ell)$. Since integrality of eigenvalues can be computed on a
finite extension of $K$, 
 we may assume that $V$ is smooth. 
If $K$ has characteristic zero, there is a good compactification
$j:V\hookrightarrow W$, with $W$ smooth proper over $K$ and $D=W\setminus V=
\cup D_i$ a
strict normal crossing divisor. Then the long exact sequence \setcounter{equation}{0}
\ga{a.5}{
\ldots \to H^i_D(\bar{W}, \Q_\ell)\to H^i(\bar{W}, \Q_\ell) \to 
H^i(\bar{V}, \Q_\ell)\to \ldots}
and Theorem \ref{apthm:int} applied to the cohomology of $W$ reduce to showing integrality for $H^i_D(\bar{W}, \Q_\ell)$. As in (3.3) of the main  article, the Mayer-Vietoris spectral sequence
\ga{a.6}{E_1^{-a+1,b}=
\oplus_{|I|=a} H^{b}_{D_I}(\bar{W}, \Q_\ell)
\Rightarrow H^{1-a+b}_D(\bar{W}, \Q_\ell),}
with $D_I=\cap_{i\in I} D_i$, 
reduces to the case where $D$ is smooth projective of codimension $\ge 1$.
Then  
 purity together with
 Theorem \ref{apthm:int}
allow us to conclude. If $K$ has equal positive characteristic, 
we apply  \cite[Th.~6.5]{apdeJ} again to find $\pi: V'\to V$ generically finite and
$j:V'\hookrightarrow W$ a good compactification. As $\pi^*: 
H^i(\bar{V}, \Q_\ell)\hookrightarrow H^i(\overline{V'}, \Q_\ell)$ is injective, we conclude as above. 
\Endproof\vskip4pt  

Corollary \ref{apcor:int} gives some flexibility as we do not assume 
that $V$ is projective. In particular, 
one can apply the same argument as in the proof of Theorem 2.1 of the main article in order to show an improved version of Theorem 1.5, (ii) there:
\begin{cor} \label{apcor:int2}
Let $V$ be a smooth scheme of finite type over $K${\rm ,} and $A\subset V$ be a 
codimension $\kappa$ subscheme. 
Then the eigenvalues of $\Phi$ on $H^i_A(\bar{V}, \Q_\ell)$ are divisible by $|k|^\kappa$ as algebraic
integers.
\end{cor}

\Proof One has a stratification $\ldots \subset A_i
\subset A_{i-i} \subset \ldots A_0=A$ by closed subschemes defined over $K$ with $A_{i-1}\setminus A_i$ smooth.  
The $\Phi$-equivariant long exact sequence
\ga{a.7}{
\ldots \to H^m_{A_i}(\bar{V}, \Q_\ell)\to H^m_{A_{i-1}}(\bar{V}, \Q_\ell)
\to  H^m_{(A_{i-1}\setminus A_i)}(\overline{V\setminus A_i}, \Q_\ell)
\to \ldots
}
together with purity and  Corollary \ref{apcor:int} allow us to conclude
 by induction on the codimension. 
\hfill\qed

\begin{rmk}
One has to pay attention to the fact that even if Theorem \ref{apthm:int}  generalizes
 Theorem 1.5 i) of the article to $V$ not necessarily smooth, there is no such generalization of  Theorem 1.5 ii) to the 
nonsmooth case, even on a finite field. Indeed, let $V$ be a  rational curve 
with one node. Then $H^1(\bar{V}, \Q_\ell)=\Q_\ell(0)$ as we see from the normalization sequence, yet $H^1_{{\rm node}}(\bar{V}, \Q_\ell)=
H^1(\bar{V}, \Q_\ell)$ as the localization map $H^1(\bar{V}, \Q_\ell)\to 
H^1(\overline{V \setminus {\rm node}}, \Q_\ell) $ 
factor through $$H^1(\overline{{\rm normalization}}, \Q_\ell)=0.$$  
So we cannot improve the integrality statement to a divisibility statement in 
general. In order to force divisibility, one needs the divisor supporting the cohomology to be in good position with respect to the singularities. 
\end{rmk}

\end{document}